\documentclass[12pt]{article}
\title{An example of a covering surface with\\ movable natural boundaries
}
\author{Claudi Meneghin
}
\begin{document}
\maketitle
\bibliographystyle{plain} 
\parindent=8pt

\def\CIRC{\mathop{\tt o}\limits}
\def\quan{\vrule height6pt width6pt depth0pt}
\def\QUAN{\ \quan}
\def\R{\hbox{\tt R}}

\newtheorem{definition}{Definition}
\newtheorem{defi}[definition]{Definition}
\newtheorem{lemme}[definition]{Lemme}
\newtheorem{proposition}[definition]{Proposition}
\newtheorem{theorem}[definition]{Theorem}        
\newtheorem{corollaire}[definition]{Corollaire}  
\newtheorem{remark}[definition]{Remark}
  
\font\sdopp=msbm10 scaled\magstep1
\font\xdopp=msbm10 
\def\ER {\sdopp {\hbox{R}}}
\def\CI {\sdopp {\hbox{C}}}
\def\CIP {\xdopp {\hbox{C}}}
\def\DI {\sdopp {\hbox{D}}}
\def\EE{\sdopp {\hbox{E}}}
\def\EN{\sdopp {\hbox{N}}}
\def\ZETA{\sdopp {\hbox{Z}}}
\def\PI {\sdopp {\hbox{P}}}
\def\M{\hbox{\tt\large M}}
\def\N{\hbox{\tt\large N}} 
\def\P{\hbox{\boldmath{}$P$\unboldmath}} 
\def\f{\hbox{\large\tt f}} 
\def\e{\hbox{\boldmath{}$e$\unboldmath}} 

\def\F{\hbox{\boldmath{}$F$\unboldmath}} 
\def\H{\hbox{\boldmath{}$H$\unboldmath}} 
\def\h{\hbox{\large\tt h}} 
\def\id{\hbox{\boldmath{}$id$\unboldmath}} 
\def\pphi{\hbox{\boldmath{}$\phi$\unboldmath}} 
\def\Ef{\phi}
\def\ppsi{\hbox{\boldmath{}$\psi$\unboldmath}} 
\def\Ppsi{\hbox{\boldmath{}$\Psi$\unboldmath}} 
\def\P{\hbox{\boldmath{}$P$\unboldmath}} 
\def\Y{\hbox{\boldmath{}$Y$\unboldmath}} 
\def\tr{\hbox{\boldmath{}$tr$\unboldmath}} 
\def\f{\hbox{\boldmath{}$f$\unboldmath}} 
\def\u{\hbox{\boldmath{}$u$\unboldmath}} 
\def\ii{i} 
\def\e{\hbox{\boldmath{}$e$\unboldmath}} 
\def\v{\hbox{\boldmath{}$v$\unboldmath}} 
\def\U{\hbox{\boldmath{}$U$\unboldmath}} 
\def\V{\hbox{\boldmath{}$V$\unboldmath}} 
\def\W{\hbox{\boldmath{}$W$\unboldmath}} 
\def\id{\hbox{\boldmath{}$id$\unboldmath}} 
\def\alph{\hbox{\boldmath{}$\aleph$\unboldmath}} 
\def\bet{\hbox{\boldmath{}$\beta$\unboldmath}} 
\def\gam{\hbox{\boldmath{}$\gamma$\unboldmath}} 
\def\U{\mathop{u}\limits}
\def\f{\hbox{\boldmath{}$f$\unboldmath}} 
\def\g{\hbox{\boldmath{}$g$\unboldmath}} 
\def\h{\hbox{\boldmath{}$h$\unboldmath}} 
\def\labelle #1{\label{#1}}
\def\obmet #1{  }

\begin{abstract}
We exhibit a covering surface of the punctured complex plane 
(with no points over $0$)
whose natural projection mapping fails to be a topological 
covering, due to the existence
of branches with natural boundaries, projecting over different 
slits in $\CI^*$.
\end{abstract}


In \cite{beardon}, A.F.Beardon shows an example
of a 'covering surface' over a region in the unit disc
which is not a topological covering (We recall that a 'covering surface' of a region $D\subset\CI$ 
is a Riemann surface $S$ admitting 
a surjective conformal mapping onto $D$, among which there are
the {\it analytical continuations} of holomorphic germs;
a continuous mapping $p$
of a topological space $Y$ onto another one
$X$ is a  'topological covering' provided that
each point $x\in X$ admits an open neighbourhood 
$ {\cal U}  $ such that the restriction of $p$
to each connected component  $ {\cal V}_i  $
of $p^{-1}({\cal U}  )$ is a homeomorphism  of $ {\cal V}_i  $ 
onto ${\cal U}$, see \cite{forster, klaus}).

Beardon's example takes origin by the analytical continuation of
a holomorphic germ, namely a branch of the inverse of 
the Blaschke product
$$
B(z)=\prod_{n=1}^{\infty}
\left(
\frac{z-a_n}{1-\overline{a}_n z}
\right),
$$
with suitable hypotheses about the $a_n$'s.
It is shown that there is a sequence of inverse branches
$f_k$ of $B$ at $0$ such that the radii of convergence of
their Taylor developments tend to $0$ as $k\to\infty$.
This of course prevents the natural projection (taking a germ at $z_0$ to the point $z_0$) of the analytical
continuation of any branch of $B^{-1}$ to be a topological covering.

Quoting Beardon, such an example shows that 
{\it there is a significant 
difference between the definition of a covering surface used by 
complex analysts and that used by topologists}.

We remark that, by another point of view, it is easy to contruct
examples of covering surfaces having some kind of isolated 
'missing points' (i.e. singular points for the analytical 
continuation), preventing their projection mappings from being 
topological coverings. 

For instance, consider the entire function 
defined by
$g(z):=\exp\exp(z) $. We can easily define infinitely 
many branches of its 'inverse' in a neighbourhood of any point
in $\CI^*$ by 
choosing suitable branches of the 'logarithm function'. 
Let $\gamma :I\to\CI^*$ be defined by
$\gamma (t)= t+(1-t)\cdot e$. 

Now let us construct the following branches of {\bf log}:
let $\ell_1$ be the branch 
taking the value $1$ at $w=e$, $\ell_2$ be the branch defined 
on $\DI(1,1)$ by extending the Taylor development of the real 
logarithm function, centered at $1$ (thus $\ell_2(1)=0$) and 
$\ell_3:= \ell_1+2\pi i $. Now both $\ell_2\circ\ell_1$ and
$\ell_2\circ\ell_3$ are inverse branches of $g$ and both admit 
analytical continuation across $\gamma \vert_{[0,1)}$.

Clearly $\ell_2\circ\ell_1$ does not admit analytical continuation up to $1=\gamma (1)$, wheras $\ell_2\circ\ell_3$ does; thus, again,
the natural projection is not a topological covering.

One could ask if we can construct such examples 
dealing with nonisolated singularities. 
In this note, we aim at showing that the answer
is yes: in fact we construct a covering surface of
$\CI^*$
(once more originating from the analytical continuation of a 
holomorphic germ),  admitting singular-point arcs 
over the slits $\EE_n:=\DI(0, e^n)\cup[e^{n-1},e^n]$,
$n\in\ZETA$.

By imagery, we could thought of these
singularities as 'movable natural boundaries' (not to be confused 
with those arising in the theory of complex ordinary differential
equations, see e.g. \cite{chazy}). These singularities will appear 
(or 'disappear') at $\EE_n$, after a quantity connected to
the winding number around $0$ of the path
which we carry out analytical continuation along (see definition
\ref{wind}). 

Of course, this prevents the natural projection of the analytical continuation

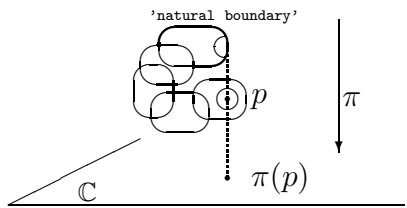
\begin{tabular}{cc}
\begin{picture}(150,48)(10,30)
\put(10,0){\line(1,0){150}}
\put(10,0){\line(2,1){50}}
\multiput(93,10)(0,2){25}{\line(0,1){1}}
\put(135,70){\vector(0,-1){50}}
\put(92,60){\oval(8,8)[l]}
\put(93,40){\circle*{2} }
\put(93,40){\circle{8} }
\put(105,40){\makebox(0,0){$p$}}
\put(93,10){\circle*{2} }
\put(113,10){\makebox(0,0){$\pi(p)$}}
\put(65,43){\oval(15,20)}
\put(75,35){\oval(22,15)}
\put(90,40){\oval(20,15)}
\put(80,60){\oval(26,15)}
{\thicklines
\put(80,60){\oval(26,15)[r]}
\put(80,60){\oval(26,15)[tl]}
}
\put(90,71){\makebox(0,0){\hbox{
\tt\tiny 'natural boundary'}}}
\put(70,53){\oval(20,15)}
\put(40,5){\makebox(0,0){$\CIP$}}
\put(140,40){\makebox(0,0){$\pi $}}
\end{picture}
&
\begin{minipage}{213pt}
from being a topological
covering (see figure aside).
To start constructing our example, let us consider the
power series
$h(z):=\sum_{\nu =0}^{\infty}z^{2^\nu }=1+z^2+z^4+
z^8+...$, defining a holomorphic function in $\DI$ and having a
natural boudary at $\partial \DI$.
\end{minipage}
\end{tabular}
\ \\
Now, let 
$
\displaystyle
{\cal E}_n:=\{x+iy\in\CI\colon n\leq
x \leq n+1\,\hbox{\rm and}\, y> 2\pi n\}\  (n\in \ZETA)
$ and
${\cal E}:=\bigcup_{n=-\infty}^{\infty}
{\cal E}_n$.
Since ${\cal E}   $ is simply connected, by Riemann 
mapping theorem,  

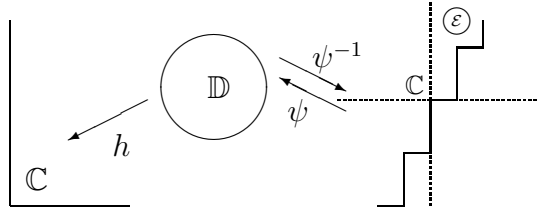
\begin{tabular}{cccc}
\begin{picture}(30,50)(0,30)
\put(0,0){\line(0,1){70}}
\put(0,0){\line(1,0){45}}
\put(10,10){\makebox(0,0){$
 \CI$}}
\end{picture}
&
\begin{picture}(70,50)(0,30)
\put(0,23){\makebox(0,0){$h$}}
\put(10,40){\vector(-2,-1){30}}
\put(35,45){\circle{40}}
\put(35,45){\makebox(0,0){
 \DI}}
\put(80,56){\makebox(0,0){
 $\psi^{-1}$}}
\put(60,56){\vector(2,-1){25}}
\put(85,36){\vector(-2,1){25}}
\put(65,35){\makebox(0,0){
 $\psi$}}
\end{picture}
&
\begin{picture}(80,50)(10,30)
\put(70,45){\makebox(0,0)
}
\multiput(45,0)(0,2){39}{\line(0,1){1}}
\multiput(10,40)(2,0){39}{\line(1,0){1}}
\put(25,0){\line(1,0){10}}
\put(35,0){\line(0,1){20}}
\put(35,20){\line(1,0){10}}
\put(45,20){\line(0,1){20}}
\put(45,40){\line(1,0){10}}
\put(55,40){\line(0,1){20}}
\put(55,60){\line(1,0){10}}
\put(65,60){\line(0,1){10}}
\put(39,46){\makebox(0,0){\tiny $\CIP$}}
{\thinlines
\put(55,70){\circle{10}}
\put(55,70){\makebox(0,0){\tiny ${\cal E}$}}
}
\end{picture}
&
\begin{minipage}{160pt}
\vskip.1cm\noindent
there exixts a biholomorphic map
$\psi\colon{\cal E}\rightarrow \DI   $, thus
$ h\circ\psi  $
has a natural boundary at 
$ \partial{\cal E}$.

Now we proceed to define a slight generalisation of
the notion of {\it winding number}, aimed
at coping with nonclosed paths:
\end{minipage}
\end{tabular}

\begin{defi}
Let $ \gamma\colon[0,1]
\rightarrow \CI\setminus\{ 0  \}  $ be a path: 
we say that $ \gamma $ {\sf winds n times around 0}
if 
$$
\left[
\frac{
\Im\left(
\ell_1\circ\gamma (1)
-\ell_0\circ\gamma (0)
\right)    
}{2\pi }    \right]=\hbox{\sf n},
$$
where $ ({\cal U}_0,\ell_0)  $ is any branch of the logarithm
in a neighbourhood of $ \gamma(1)  $,
$ ({\cal U}_1,\ell_1)  $ the branch got by analytical
continuation across $ \gamma  $ and \boldmath $ [\ ]  $
denotes the {\rm integer-part} operator. We shall write 
$\hbox{\sf W}(\gamma)=\hbox{\sf n}$.
\label{wind}
\end{defi}

In the following theorem, the main construction of this note
will be fully depicted: let $\lambda $ be the branch of the complex 
'logarithm function' taking
the value $\hbox{ln}(1/2)$ at $1/2$ (defined by its 
Taylor development in $\DI(1/2,1/2)$) and 
$f :=h\circ\psi\circ\lambda$. Let $\gamma $ be a path as in 
definition \ref{wind}.
We have:

\begin{theorem}
{\tt (a) - conditioned everywhere continuability}: for every
$\omega \in \CI^*$ there exists a path $\beta $ from $1/2$ to
$\omega$ along which $f$ can be continued;
{\tt (b) - movable boundaries:} let $M,N\in\ZETA$; for every 
$\omega\in\partial \DI(0,e^M)\cup [e^{M-1},e^M] $ let 
$\varphi $ be a path joining $1/2$ and $\omega$ in $\CI^*$ such
that $\hbox{\sf W}(\gamma)=\hbox{\sf N}$ and $f$ admits analtyical
continuation along $\varphi\vert_{[0,1)}$. If $N>M$, then $f$
admits analytical continuation along $\varphi $ up to
$\omega =\varphi (1)$; if $N=M$, then $f$ does not.
\labelle{akill}
\end{theorem}
{\bf Proof}\\ 
{\tt (a)} let $\omega \in \CI^*$:
since the exponential function
admits the period $2\pi i$, we can find $\zeta\in {\cal E} $
such that $e^{\zeta }=\omega $. 
Let $\alpha $ be a curve in ${\cal E} $ joining 
$\hbox{ln}(1/2)$ and $\zeta$; the thesis 
follows now by setting $\beta := \exp(\alpha )$:
$\lambda $ can of course be continued along $\beta$ in such a way 
to get a finite chain of holomorphic function elements 
$({\cal U}_i, \Lambda_i)_{i=0...M}$ such that 
$\bigcup_{i=0}^N\Lambda_i({\cal U}_i)\subset {\cal E}$; the analytical continuation
of $f$ along $\beta$ now follows by composing the $\Lambda_i$'s
with $h\circ\psi $ on the left.
\ \\ 
{\tt (b)}: Now, $\lambda$ can be certainly continued 
across $\varphi$; let $\ell $ be the 
branch of the logarithm
obtained in this way at $\varphi (1)$. 
Since $f$ admits analytical continuation along 
$\varphi\vert_{[0,1)}$, we have 
$\ell\circ\varphi\vert_{[0,1)}\subset {\cal E}$. 
Now $\hbox{\sf W}(\varphi)= N$, so if 
$\omega \in \partial \DI(0,e^M) $, then $\Re(\ell(\omega))=M$
and $2\pi (N-1)\leq \Im(\ell(\omega))< 2\pi N $; 
if $\omega \in [e^{M-1},e^M]$, then 
$M-1\leq \Re(\ell(\omega))\leq M$ and 
$\Im (\ell(\omega))=2\pi (N-1)$.
Hence, if $N>M$, then $\ell(\omega)\in{\cal E}$, implying that 
$h\circ\psi\circ\ell$ admits analytical continuation up to
$\omega =\varphi (1) $; by consequence, so does $f$; 
if $N=M$, then $\ell(\omega)\in\partial {\cal E}$,
therefore we cannot carry out analytical continuation
of $h\circ\psi$ at $\ell\circ \varphi(1)$, which
means that $f$ could not be continued at $\omega$.
\QUAN
\vskip.5cm
The picture in theorem \ref{akill} is clear: there are some kind
of 'natural boundaries' for the analytical continuation of $f$
which can be 'pushed farther' by winding the 
analytical-continuation path a suitable number of times around 
$0$; as already stated, 
this prevents the natural projection of the analytical 
continuation of $f$ from being a topological covering.


\begin{thebibliography}{9}

%
\bibitem{beardon}
A.F.Beardon
{\it A remark on analytic continuation}
{Proceedings of the AMS, 128, 5}
%
%
%
\bibitem{forster}
{Otto Forster}
{\it Lectures on Riemann surfaces}
{Springer Verlag, 1981}
%
\bibitem{klaus}
{Klaus J\"anich}
{\it Topology}
{Springer Verlag, 1994}
%
%

\end{thebibliography}
\end{document}